\documentclass[12pt]{article}
 \usepackage{graphicx}
 \usepackage{amssymb}
 \usepackage{epstopdf}
 \newtheorem{theorem}{Theorem}
 \newtheorem{lemma}[theorem]{Lemma}
 \newtheorem{corollary}[theorem]{Corollary}
 \newtheorem{proposition}[theorem]{Proposition}
 
 \newtheorem{question}[theorem]{Question}
 \def \Z{\mathbb Z}
 
 \def \HH {\mathrm{H}}
 \def \C {\mathbb C}
 \def\Q{\mathbb Q}
 
 \def \R{\mathbb R}
 
 \def \sS{\cal S} 
 
 \def\sK{\cal K}
 \def\sL{\cal L}
 
 \def\link{\ell \mathrm{k}}
 \def\qed{\hfill\framebox(5,5){}\newline\medskip}
 \def\sign{\mathrm {s}}
 \def\tr{\mathrm{tr}}
 \def\iso{\simeq}

 \title{Salem-Boyd sequences and Hopf plumbing}
 \author{Eriko Hironaka\thanks
{The author was supported in part by a JSPS Research Fellowship at Osaka University}
} 
\begin{document}
  \maketitle
 \abstract {Given a fibered link, consider the characteristic polynomial of the monodromy
 restricted to first homology.  This generalizes the notion of the
 Alexander polynomial of a knot.  We define a construction, called iterated plumbing, to create 
 a sequence of fibered links from a given one.
The resulting sequence of characteristic polynomials for these links has the same
 form as those arising in work of Salem and Boyd in their study of distributions of 
 Salem and P-V numbers.    From this we deduce information about the asymptotic
 behavior of the large roots of the generalized Alexander polynomials, 
 and define a new poset structure for Salem fibered links.}

 \section{Introduction}  

 Let $(K,\Sigma)$ denote a fibered link $K \subset S^3$ with fibering surface $\Sigma$.
 Hopf plumbing defines a natural operation on fibered links that allows one to construct 
new fibered links from a given one while keeping track of useful information \cite{Stallings:fibered} \cite{Gabai:Murasugi}.  Furthermore,
a theorem of Giroux \cite{Giroux02} shows that any fibered link can be obtained from
the unknot by a sequence of Hopf plumbings and de-plumbings (see also \cite{Harer:fibered}).

A fibered link $(K,\Sigma)$ has an associated
homeomorphism $h : \Sigma \rightarrow \Sigma$, called the {\it monodromy}
of $(K,\Sigma)$, such that the complement in $S^3$ of a regular
neighborhood of $K$ is homeomorphic to a mapping torus for $h$.    Let $h_*$
  be the restriction of $h$ to first homology $\HH_1(\Sigma,\R)$, and let
  $\Delta_{(K,\Sigma)}(t)$ be the characteristic polynomial of the monodromy $h_*$.  
  If $K$ is connected, that is, a {\it fibered knot}, then $\Delta_{(K,\Sigma)}(t)$ is 
  the usual Alexander polynomial of $K$ and the mapping torus structure is unique.
 We extend this terminology and call 
  $\Delta_{(K,\Sigma)}(t)$ the {\it Alexander polynomial} of the fibered link $(K,\Sigma)$. 
      
  A polynomial $f$ of degree $d$
  is {\it reciprocal} if  $f=f_*$, where $f_*(t) = t^df(1/t)$.  The Alexander polynomials $\Delta_{(K,\Sigma)}(t)$
 are monic integer polynomials and reciprocal up to multiples of $(t-1)$.   
Burde  \cite{Burde:Alex} shows that there exists a fibered knot $(K,\Sigma)$
with $\Delta_{(K,\Sigma)} = f$, if and only if
 \begin{description}
 \item{(i)} $f$ is a reciprocal monic integer polynomial; and
 \item{(ii)} $f(1) = \pm 1$,
 \end{description}
Kanenobu \cite{Kanenobu81} shows that (i) is true if and only if
$\Delta_{(K,\Sigma)} = f$
up to multiples of  $(t-1)$, where $(K,\Sigma)$ is a fibered link.
 Our goal in this paper is to study how the roots 
  of $\Delta_{(K,\Sigma)}(t)$ are affected by Hopf plumbing.
 
 In Section~\ref{plumbing-section}, we define a construction called iterated (trefoil) plumbing,
 which produces a sequence of fibered links $(K_n,\Sigma_n)$ from a given fibered link 
 $(K,\Sigma)$ and a choice of path $\tau$ properly embedded on $\Sigma$, 
 called the {\it plumbing locus}.
 
 Our main result is the following.
 
 \begin{theorem}\label{main-theorem}
  If $(K_n,\Sigma_n)$ is obtained from $(K,\Sigma)$ by $\pm$ iterated trefoil plumbing,
 then there is a polynomial $P = P_{\Sigma,\tau}$ depending only on the location and
orientation of the plumbing, such that $\Delta_n = \Delta_{(K_n,\Sigma_n)}$ is given by \begin{eqnarray}\label{Salem-Boyd-eqn}
 \Delta_n(t) = \frac{t^{2n} P(t) \pm (-1)^{r} P_*(t)}{t+1},
 \end{eqnarray}
 where $r$ is the number of components of $K$.
  \end{theorem}
 
We  call sequences of polynomials of the form given in 
 Equation~\ref{Salem-Boyd-eqn} {\it Salem-Boyd sequences}, after work of Salem \cite{Salem44}
 and Boyd \cite{Boyd77} \cite{Boyd89}.  
  
 For a monic integer polynomial $f(t)$, let $\lambda(f)$ be the maximum absolute value
 among all roots of $f(t)$; $N(f)$, the number of roots with absolute value greater than one;
 and $M(f)$, the product of absolute values of roots of $f$ whose absolute value is greater
 than one.   The latter invariant $M(f)$ is known as the {\it Mahler measure of $f$}.
 Clearly $N(f)$ is discrete, while $\lambda(f)$ can be made arbitrarily close to but greater than
 one, for example, by taking $f(t) = t^n - 2$.    Whether or not the values of $M(f)$ can
 also be brought arbitrarily close to one from above is an open problem posed
 by Lehmer in 1933 \cite{Lehmer33}.  Lehmer originally formulated his question as follows:
 
 \begin{question}[Lehmer] For each $\delta > 0$ does there exist a monic integer polynomial
 $f$ such that $1 < M(f) < 1+ \delta$?
 \end{question}

 We are still far from answering Lehmer's question, but  show in 
 Section~\ref{Salem-Boyd-section} how to apply Salem and Boyd's work and
  Theorem~\ref{main-theorem}
 to obtain information about the asymptotic behavior of $N(\Delta_n)$,  
 $\lambda(\Delta_n)$, and $M(\Delta_n)$ from properties of the original 
 fibered link and location of plumbing.    
 
 \begin{theorem}\label{convergence-theorem} The sequences $N(\Delta_n)$, 
 $\lambda(\Delta_n)$ and  $M(\Delta_n)$ converge to $N(P)$, $\lambda(P)$,
 and $M(P)$, respectively, where $P = P_{\Sigma,\tau}$.
 \end{theorem}
 
Theorem~\ref{convergence-theorem} is useful for finding 
minimal Mahler measures appearing in particular families of fibered links, since the polynomials 
 $P_{\Sigma,\tau}$ are easy to compute for explicit examples.  We give an illustration
 in Section~\ref{example-section}.

 Iterated plumbing may be seen as the result of iterating full twists  on a pair of strands
 of $K$, with some extra conditions on the pair of strands.  For the case where $K$ has
 one component,  the  convergence of Mahler
 measure in Theorem~\ref{convergence-theorem} agrees with a result of 
 Silver and Williams, which in general form may be stated as follows.  Let $L$
 be a link and $k$ an unknot disjoint from $L$ such that $L$ and $k$ have non-zero
 linking number.  Let $L_n$ be obtained from $L$ by doing $1/n$ surgery along $k$.
 This amounts to taking the strands of $L$ encircled by $k$ and doing $n$ full-twists
 to obtain $L_n$.
 Silver and Williams show that
 the multi-variable Mahler measures of the links $L_n$ 
 converge  to the multi-variable Mahler measure of $L \cup k$
 \cite{S-W:Mahler}.    Combining our results with that of Silver and Williams,
 and using the formulas for $P_{\Sigma,\tau}$ given in 
 Section~\ref{plumbing-section} (Equations~\ref{FirstP-equation} and \ref{SecondP-equation})
gives a new effective way to calculate the multi-variable Mahler measure of $L \cup k$.
   
 It is not hard to see that if one fixes
 the degree of $f$, then the answer to Lehmer's question is negative.  
 Theorem~\ref{convergence-theorem} makes it possible to study Mahler measures
 for sequences of fibered links
 whose fibers have increasing genera, and hence for polynomials of increasing degree.
Although, in general, $\lambda(\Delta_n)$ and $M(\Delta_n)$ are not monotone sequences
(see Theorem~\ref{Arg-Theorem}),
monotonicity can be shown (at least for large enough $n$) when $P_{\Sigma,\tau}$ has 
special properties.

In Section~\ref{Salem-Boyd-section}, we review properties of Salem-Boyd sequences,
following work of Salem \cite{Salem44} and Boyd \cite{Boyd77}, and consider the question
of monotonicity. 
A {\it Perron polynomial} is a monic integer polynomial $f$
with a real root $\lambda = \lambda(f) > 1$ satisfying $|\alpha| < \lambda$ for all roots $\alpha$
of $f$ not equal to $\lambda$.  

 \begin{theorem}\label{PerronCase-Theorem}  
 Suppose $P_{\Sigma,\tau}$ is a Perron polynomial.  Then $\lambda(\Delta_n)$
 is an eventually monotone (increasing or decreasing)
 sequence converging to $\lambda(P_{\Sigma,\tau})$.
 \end{theorem}
 
 \noindent
 In the special case when $N(P_{\Sigma,\tau}) = 1$, more can be shown by applying results of
 Salem \cite{Salem44} and Boyd \cite{Boyd77}.   
    
 \begin{theorem}\label{SalemCase-Theorem}   Suppose $N(P_{\Sigma,\tau}) = 1$.  
 Then $M(\Delta_n) = \lambda(\Delta_n)$ is a
  monotone (increasing or decreasing) sequence converging to $\lambda(P_{\Sigma,\tau})$.
 \end{theorem}

 Section~\ref{applications-section} studies the poset structure on fibered links defined by Hopf
 plumbing, and the corresponding poset structure on homological dilatations.  We also give
an example in Section~\ref{applications-section} that shows how
 Theorem~\ref{SalemCase-Theorem} can be used to give explicit solutions to Lehmer's
 problem for restricted families.
  
 \medskip
 \noindent
 {\sc Acknowledgements:}  I am indebted to J.S.P.S. who funded my research, and the staff of the
 Osaka University Mathematics Department and 
 my host Makoto Sakuma  for their kind hospitality and support during the writing 
  of this paper.  I would also like to thank  
 D. Silver and S. Williams for many helpful discussions, 
 and S. Williams and E. Kin for bringing my attention to the example in Section~\ref{example-section}.
 
\section{Iterations of Hopf plumbings.}\label{plumbing-section}

We recall some basic definitions surrounding the Alexander polynomial of
an oriented link.  A {\it Seifert surface} for an oriented link $K$ is an oriented surface 
$\Sigma$ whose boundary
is $K$.   For any collection of free loops $\sigma_1,\dots,\sigma_d$ on $\Sigma$ forming a basis for
 $\HH_1(\Sigma;\R)$, the  associated
{\it Seifert matrix} $S$ is given by 
$$
S = [\link (\sigma_i^+,\sigma_j)], 
$$
where $\sigma_i^+$
is the push-off of $\sigma_i$ off $\Sigma$ into $S^3 \setminus \Sigma$ in the positive 
direction with respect to the orientation of $\Sigma$, and $\link(,)$ is the
linking form on $S^3$.  Let $S^\tr$ denote the transpose of $S$.   The
polynomial
$$
\Delta_K(t) = |tS - S^\tr|
$$
is uniquely defined up to units in the Laurent polynomial ring $\Lambda(t) = \Z[t,t^{-1}]$, 
and is reciprocal (it is the same if
$S$ is replaced by $S^\tr$).   For the purposes of this paper, we will always normalize 
$\Delta_K$ so that 
$\Delta_K(0) \neq 0$, and the highest degree coefficient of $\Delta_K(t)$ is positive.
Then for any nonsingular Seifert matrix for $K$, 
$$
\Delta_K(t) = \sign(S) |tS - S^\tr|,
$$
where $\sign(S)$ is the sign of the
coefficient of $|tS-S^\tr|$ of highest degree. 

If $K$ is fibered, and $\Sigma$ is the fibering surface, then the Seifert matrix $S$ is
invertible over the integers, and the monodromy
restricted to $\HH_1(\Sigma;\R)$ satisfies $h_* = S^\tr S^{-1}$.
In this case $\sign(S) = |S|$, and $\Delta_K(t)$ is characteristic polynomial of $h_*$.      
Since $|S|$ is invariant under change of basis, and the fiber surface is fixed, we will
write $\sign(K) = \sign(S)$ if $K$ is fibered.
If $K$ is a fibered knot, then $\sign(K) = \Delta_K(1)$.
\begin{figure}[htbp]
\begin{center}
\includegraphics[width=3in]{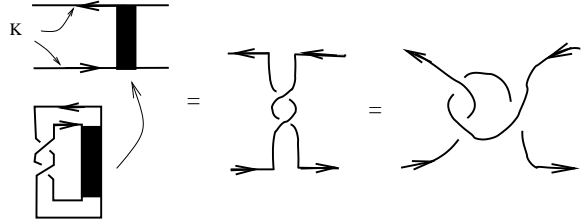}
\caption{Positive Hopf plumbing}
\label{Hopfplumbing-fig}
\end{center}
\end{figure}

A {\it properly embedded path} on $\Sigma$ is a smooth embedding
$$
\tau: [0,1] \rightarrow \Sigma
$$
such that $\tau(0),\tau(1) \in \partial \Sigma$. 
The surface $\Sigma^+_2 (\tau)$ (resp., $\Sigma^-_2$ is obtained from $\Sigma$ 
by positive (resp., negative) 
Hopf plumbing if it is obtained from $\Sigma$ by gluing on a positive (resp., negative)
Hopf band as in Figure~\ref{Hopfplumbing-fig}.  The definition is independent of
the orientation of $\tau$.

\begin{figure}[htbp]
\begin{center}
\includegraphics[height=0.5in]{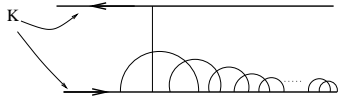}
\caption{Base paths for iterated Hopf plumbings}
\label{basepaths-fig}
\end{center}
\end{figure}

Set  $\Sigma^\pm_1 = \Sigma$.  
For $n \ge 1$, let $\Sigma^\pm_{n+1}$
be the  {\it (positive or negative) 
Hopf $n$-plumbing} of $\Sigma$ 
along  $\tau$, which is obtained by
Hopf plumbing along $n$ paths as shown in Figure~\ref{basepaths-fig}, starting with
the vertical path $\tau$.

\begin{figure}[htbp]
\begin{center}
\includegraphics[height=1.25in]{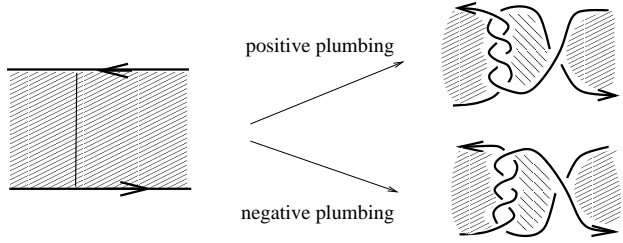}
\caption{Result of iterated Hopf $4$-plumbing}
\label{iteratedplumbing-fig}
\end{center}
\end{figure}

The positive (resp., negative) Hopf $n$-plumbings can also be considered as a Murasugi sum
of $\Sigma$ with the fiber surface of the torus link $T(2,n)$ (resp., $T(2,-n)$).  
Let $K^\pm_n(\Sigma,\tau)$ be the boundary of the surface $\Sigma^\pm_n$.  
For $n=1$, we have 
$K^\pm_1=K$.  The local oriented link diagram  for $K^\pm_n$
is shown in Figure~\ref{iteratedplumbing-fig}, and $\Sigma^\pm_n$ is the corresponding 
 Seifert  surface.

Denote by $\langle,\rangle$ the intersection pairing
$$
\HH_1(\Sigma;\Z) \times \HH_1(\Sigma,\partial\Sigma;\Z) \rightarrow \Z,
$$
and let $v  \in \HH_1(\Sigma;\Z)$ be the vector such that $v^\tr$ represents
the vector $[\tau]$ in $\HH_1(\Sigma;\partial\Sigma;\Q) \iso \HH_1(\Sigma;\Q)^{\mathrm{dual}}$.  
Then, in terms of the basis $\sigma_1,\dots,\sigma_d$, $v$ is given by
$$
v = [\langle \sigma_1,\tau\rangle, \dots, \langle \sigma_d,\tau \rangle].
$$
Set
\begin{eqnarray}\label{FirstP-equation}
P^\pm_{\Sigma,\tau}(t) =  |tI - (S^\tr \mp  vv^\tr)S^{-1}|,
\end{eqnarray}
where $I$ is the identity matrix.

Before proving the Theorem~\ref{main-theorem}, we put $P^\pm_{\Sigma,\tau}$ in an alternate form.
Let $N(\tau)$ be a regular neighborhood of $\tau$ on $\Sigma$, and let 
$\Sigma_0 = \Sigma \setminus N(\tau)$.  Let $K_0 = K_0(\Sigma,\tau) = \partial \Sigma_0$.
Let $\sigma_1,\dots,\sigma_{d-1}$ be a collection of free loops on $\Sigma_0$
forming a basis for $\HH_1(\Sigma_0;\Z)$.  Let $\sigma_d$ be a free loop on
$\Sigma$ so that $\sigma_1,\dots,\sigma_d$ is a basis for 
$\HH_1(\Sigma;\Z)$, and such that
$\langle \sigma_d,\tau \rangle = 1$.    
Let $S_1$ and $S_0$ be the corresponding Seifert matrices for $K$ and $K_0$,
respectively.      

\begin{lemma}\label{nonsingular-lemma} The Seifert matrix $S_0$ is non-singular.
\end{lemma}

\noindent{Proof.}  By our definitions, the transpose of the Seifert matrix defines a linear transformation 
from the first homology of the Seifert surface to its dual.  We thus have a commutative diagram
$$
\begin{array}{ccc}
\HH_1(\Sigma_0;\Z) &{\shortstack {$S_0^\tr$ \\ $\rightarrow$ }}& \HH_1(\Sigma_0;\Z)^{\mathrm{dual}}\\
\vector(0,-1){12} & &\vector(0,-1){12}\\
\HH_1(\Sigma_1;\Z) &{\shortstack {$S_1^\tr$ \\ $\rightarrow$}} &\HH_1(\Sigma_1;\Z)^{\mathrm{dual}}.
\end{array}
$$
where vertical arrows are the inclusions determined by the choice of bases.
Since $S_1$ is non-singular, it follows that $S_0$ must also be non-singular.\qed

\begin{lemma}\label{propertyofP-lemma} The polynomial in Equation~\ref{FirstP-equation}
can be rewritten as
\begin{eqnarray}\label{SecondP-equation}
P^\pm_{\Sigma,\tau}(t)
=\Delta_K(t) \pm \sign(K)\sign(S_0)\Delta_{K_0}(t).
\end{eqnarray}
\end{lemma}

\noindent{Proof.}  The choice of basis $\sigma_1,\dots,\sigma_d$ above
yields the Seifert matrix 
$$
S_1 = \left [
\begin{array}{c|c}
S_0 & x \\
\hline
y^\tr & s
\end{array}
\right ]
$$
for $K$,
where $x,y \in \Z^{d-1}$, and $s \in \Z$.  The vector $v$ written with respect to the dual
elements
of $\sigma_1,\dots,\sigma_d$ is given by $v = [0,\dots,0,1]^\tr$.
We thus have
$$
|tS_1 - (S_1^\tr \mp vv^\tr)| =
\left |
\begin{array}{c|c}
tS_0 - S_0^\tr & tx - y\\
\hline
ty^\tr - x^\tr & s(t-1)\pm 1
\end{array}
\right |.
$$
Therefore
$$
P^\pm_{\Sigma,\tau}(t) = \sign(K)( |tS_1-S_1^\tr| \pm |tS_0 - S_0^\tr|)
$$
and the claim follows.  \qed

For a polynomial $g$, define 
$$\overline {g}(t) = t^{-m}g(t),
$$
where $m$ is the largest power of $t$ dividing $g$.
Then it is easy to check that  $g_*(t) = \overline{g}_*(t)$.  Also,
if $g$ and $f$ are polynomials of degrees $d'$ and $d$, respectively, then 
for $h(t) = g(t) \pm f(t)$, we have
$$
h_*(t) = g_*(t) \pm t^{d'-d} f_*(t).
$$

\begin{lemma}\label{propertyofP*-lemma} Let $r$ be the number of components
of $K$, and  $P(t) = P^\pm_{\Sigma,\tau}(t)$. Then
$$
P_*(t) = (-1)^{r+1}\left (\Delta_K(t) \mp \sign(K)\sign(S_0)t\Delta_{K_0}(t)\right ).
$$
\end{lemma}
\medskip

\noindent
{Proof.}   If $d$ is the rank of $\HH_1(\Sigma;\R)$,
we have 
$$
P_*(t) = t^d\Delta_K(1/t) \pm  \sign (K)  t^d (|(1/t)S_0-S_0^\tr|).
$$
The  Alexander polynomial of a link is 
reciprocal (anti-reciprocal) if the number of components
is odd (even).  Thus, the first summand equals 
$(-1)^{r+1}\Delta_K(t)$.  Since, by Lemma~\ref{nonsingular-lemma},  $S_0$ is a non-singular matrix, 
$t$ does not divide $|tS_0-S_0^\tr|$.
It is also not difficult to check that the number of components of 
$K_0$ and $K_1$ have
opposite parity,  and the degree of $|tS_0 - S_0^\tr|$ is one less than
the degree of $\Delta_K(t)$.   We thus have
$$
t^d|(1/t)S_0 - S_0^\tr| = (-1)^r t |tS_0 - S_0^\tr| = (-1)^r \sign(S_0)t \Delta_{K_0}(t).
$$
 \qed

Theorem~\ref{main-theorem} is implied by the following stronger version.

\begin{theorem}~\label{stronger-main-theorem}  
Let $(K_n,\Sigma_n)$ be obtained by $\pm$ iterated Hopf plumbing
on a fibered link $(K,\Sigma)$ with $r$-components.  Let $\Delta_n = \Delta_{(K_n,\Sigma_n)}$,
and let $P = P^\pm_{\Sigma,\tau}$.  Then 
$$
\Delta_n(t) = \frac{t^nP(t) \pm (-1)^{r+n}P_*(t)}{t+1}.
$$
\end{theorem}
\medskip

\noindent
{Proof.}   By Lemma~\ref{propertyofP*-lemma}, we
have 
\begin{eqnarray*}
tP(t) + (-1)^{r+1}P_*(t) 
&=& t\Delta_K(t)\pm\sign(K)\sign(S_0)t\Delta_{K_0}(t)\\
&&\qquad\qquad+ (\Delta_K(t) \mp \sign(K)\sign(S_0)t\Delta_{K_0}(t))\\
&=& (t+1)\Delta_K(t).
\end{eqnarray*} 

For $m \ge 1$,  the Seifert matrix for $S^\pm_m$ is given by
$$
S^\pm_m=\left [
\begin{array}{c|c}
S^\pm_{m-1} & 0\\
\hline
w & \pm 1
\end{array}
\right ],
$$
where $w = [0,\dots,0,-1]$.
Thus, the Alexander polynomial for $K^\pm_m$ is given by
$$
\Delta_{K^\pm_m}(t) = \sign(K^\pm_m)\left |
\begin{array}{c|c}
tS^\pm_{m-1} - (S^\pm_{m-1})^\tr & -w^\tr\\
\hline
tw & \pm (t-1)\\
\end{array} \right |.
$$
It follows that for $n\ge 2$, $\Delta_{K^\pm_n}(t)$ satisfies
$$
(t+1)\Delta_{K^\pm_n}(t) =\sign(K^\pm_n)(t+1)
\left[\pm(t-1)|tS^\pm_{n-1}-(S^\pm_{n-1})^\tr| + t|tS^\pm_{n-2}-(S^\pm_{n-2})^\tr|\right ],
$$
and $\sign(K^\pm_n)=\pm\sign(K^\pm_{n-1})$.
For $n=2$, using $\sign(K) = \sign(K_1) = \pm\sign(K^\pm_2)$, we have
\begin{eqnarray*}
(t+1)\Delta_{K^\pm_2}(t)&=& 
\sign(K^\pm_2) \left[\pm (t^2-1)|tS_1 - S_1^\tr| + (t^2+t)|tS_0-S_0^\tr| \right ]\\
&=&\pm\sign(K^\pm_2)((t^2-1)|tS_1 - S_1^\tr| \pm (t^2+t)|tS_0-S_0^\tr|)\\
&=&\sign(K)t^2(|tS_1-S_1^\tr| \pm |tS_0-S_0^\tr|)\\
&&\qquad\qquad- \sign(K)(|tS_1-S_1^\tr| \mp t|tS_0-S_0^\tr|)\\
&=&t^2 P(t) + (-1)^r P_*(t)\\
&=&t^2 P(t) + (-1)^{r+2}P_*(t).
\end{eqnarray*}
If $n > 2$, we use induction, to obtain
\begin{eqnarray*}
(t+1)\Delta_{K^\pm_n}(t)
&=&\pm\sign(K^\pm_n)[(t^2-1) |tS^\pm_{n-1} - (S^\pm_{n-1})^\tr| \pm t(t+1)|tS^\pm_{n-2} - (S^\pm_{n-2})^\tr|]\\
&=&\sign(K^\pm_{n-1})[\sign(K^\pm_{n-1})(t+1)(t-1)\Delta_{K^\pm_{n-1}}(t) \\
&&\qquad\pm
\sign(K^\pm_{n-2})t(t+1)\Delta_{K^\pm_{n-2}}(t)]\\
&=&(t-1)(t+1)\Delta_{K^\pm_{n-1}}(t) + t(t+1)\Delta_{K^\pm_{n-2}}(t)\\
&=& (t-1)(t^{n-1}P(t) + (-1)^{n-1 + r}P_*(t)) + t(t^{n-2}P(t) + (-1)^{n-2+r}P_*(t))\\
&=& t^nP(t) - t^{n-1}P(t) + (-1)^{n+r-1}tP_*(t) + (-1)^{n+r}P_*(t)\\
&&\qquad + t^{n-1}P(t) + (-1)^{n+r-2}tP_*(t)\\
&=& t^nP(t) + (-1)^{n+r}P_*(t)
\end{eqnarray*}
\qed

\section{Properties of Salem-Boyd sequences.}\label{Salem-Boyd-section}

In this section we review some general properties of roots of polynomials in Salem-Boyd
sequences (see also, \cite{Salem44}, \cite{Boyd77}), and
apply them to the Alexander polynomials of iterated plumbings.

\subsection{Asymptotic behavior of roots of Salem-Boyd sequences.}

Given a monic integer polynomial $P(t)$  define
\begin{eqnarray}\label{SalemBoyd-equation}
Q^\pm_n(t) &=& t^nP(t) \pm P_*(t).
\end{eqnarray}
We will call the sequence of polynomials given in 
Equation~\ref{SalemBoyd-equation} the {\it Salem-Boyd sequence} associated to $P$.
For all positive integers $n$,
 $Q^\pm_n(t)$ is equal to a reciprocal polynomial up to a multiple of  $t-1$.  We are
 interested in the asymptotic behavior of roots of $Q^\pm_n(t)$.
 
S. Williams suggested the use of Rouch\'e's theorem to prove the following.

\begin{lemma}\label{limitofroots-lemma}
Let $P$ be a monic integer polynomial, and let  $R(t)$ be any integer polynomial,
and
$$
Q_n(t) = t^n P(t) \pm R(t).
$$
Then the roots of $Q_n(t)$ outside $C$ converge to those of $P(t)$ counting multiplicity
as $n$ increases. 
\end{lemma}

\noindent{Proof.}  Consider the rational function
$$
S_n(t) = \frac{Q_n(t)}{t^n} = P(t) \pm \frac{R(t)}{t^n}.
$$
Let $\alpha$ be a root of $P(t)$ (counted with multiplicity), and let
$D_\alpha$ be any small disk around $\alpha$ that is also strictly outside $C$
and that contains no roots of $P(t)$ other than $\alpha$.
Then $P(t)$
has a lower bound on the boundary $\partial D_\alpha$, and thus 
there exists an $n_\alpha$ depending
on $\alpha$ and $D_\alpha$ such that
$$
\left |\frac{R(t)}{t^{n}} \right | < |P(t)|
$$
on $\partial D_\alpha$ for all $n > n_\alpha$.
By Rouch\'e's theorem, it follows that for $n>n_\alpha$, 
$P(t)$ and $S_n(t)$ (and hence also $Q_n(t)$) have
$m$ roots in $D_\alpha$ counted with multiplicity.  Since the disks could be made
arbitrarily small, and there are only a finite number of roots, the claim follows.
\qed

\begin{lemma}\label{numberofroots-lemma} 
Let $P$ be a monic integer polynomial and let $Q_n(t)$
be the associated Salem-Boyd sequence.
Then $N(Q_n) \leq N(P)$ for all $n$.
\end{lemma}

A proof of this Lemma is contained in \cite{Boyd77} (p. 317), but we include it here for the
convenience of the reader.

\medskip
\noindent{Proof.}  We first assume that $P(t)$ has no roots on the unit 
circle.  This does not change the statement's generality.
To
study the roots of $Q_n(t)$ it suffices to consider the case when $P(t)$
has no reciprocal or anti-reciprocal factors, since such factors will be factors
of $Q_n$ for all $n$. 
 If $P(t)$ has a root on the unit circle, then the 
minimal polynomial of that root would be necessarily reciprocal or
anti-reciprocal, and we can factor the minimal polynomial out of $P$ and
the $Q_n$.

Consider the two variable polynomial
\begin{eqnarray}\label{deformation-equation}
Q_n(z,u) = z^nP(z) \pm uP_*(z)
\end{eqnarray}
where $z$ is any complex number and $u \in [0,1]$.  

Suppose $P(t)$ has roots
$\theta_1,\dots,\theta_s$ outside the unit circle $C$ counted with multiplicity.
Then $Q^\pm_n(z,u)$ defines an algebraic curve $z = Z(u)$ with 
branches $z_1(u),\dots,z_s(u)$ satisfying $z_i(0) = \theta_i$.
For  $z\in C$ we have $|P(z)| = |P_*(z)|$.  Now suppose
that $0<u<1$ and $1 = |z_i(u)|$.  Then
$$
1 = |z_i(u)|^n = \frac{u|P_*(z_i(u))|}{|P(z_i(u))|} = u
$$
yielding a contradiction.  Thus, by continuity
$$
|z_i(u)| > 1, 
$$
for $u \in [0,1)$.  It follows that  $Q^\pm_n(t)$ has at most $s$ roots outside $C$.\qed

Summarizing the contents of Lemma~\ref{limitofroots-lemma} and Lemma~\ref{numberofroots-lemma}
we have the following.

\begin{theorem}\label{limitofroots-theorem}  Let $P$ be a monic integer polynomial,
and let 
$$Q_n(t) = t^n P(t) \pm P_*(t).$$
 Then
\begin{eqnarray*}
N(Q_n) &\leq& N(P);\\
\lim_{n\rightarrow \infty} \lambda(Q_n) &=& \lambda(P);\quad\mbox{and} \\
\lim_{n\rightarrow \infty} M(Q_n) &=& M(P).
\end{eqnarray*}
\end{theorem}

Theorem~\ref{main-theorem} and Theorem~\ref{limitofroots-theorem} imply 
Theorem~\ref{convergence-theorem}.

A natural question is whether $M(Q_n)$ is a monotone sequence, perhaps on
arithmetic progressions, when
$P$ has more than one root outside $C$.  The proof of Lemma~\ref{limitofroots-lemma},
does not restrict the directions by which the roots of $Q_n$ outside $C$ approach those of $P$.  If
a root $\theta$ of $P$ is not real, then the root(s) of $Q_n$ approaching $\theta$ typically
rotate around $\theta$ as they converge.   More precisely, we have the following.
For $z$ a complex number, let $A = \mathrm{Arg}(z)$ be such that
$z = |z| e^{2\pi i A}$.

\begin{theorem}\label{Arg-Theorem} Let $\alpha_1,\dots,\alpha_s$ be the roots
of $P$ outside $C$.   Take $N_0$, so that $Q_n$ has $s$ roots outside
$C$ for $n \ge N_0$.  Label these roots $\alpha_i^{(n)}$, for $i=1,\dots,s$, so that
$$
\lim_{n \rightarrow \infty} \alpha_i^{(n)} = \alpha_i.
$$
  Then, there is a constant
$c$ such that for any $\delta > 0$,
and $n > N_\delta > N_0$,
$$
\mathrm{Arg}(\alpha_i^{(n)} - \alpha_i) = c + n \mathrm{Arg}(\alpha_i)  + \delta_n,
$$
where the error term $\delta_n$ satisfies $|\delta_n| < \delta$.
\end{theorem}

\noindent
{Proof.}  Let $P_1(x)$ be the largest degree monic integer factor of $P(x)$ with no roots
outside $C$.  For $i=1,\dots,s$, we have
$$
\alpha_i^{(n)} - \alpha_i = \left (\frac{1}{\alpha_i^{(n)}} \right )^n R_n,
$$
where 
$$
R_n = 
\frac{P_*(\alpha_i^{(n)})}{P_1(\alpha_i^{(n)})
(\alpha_i^{(n)}-\alpha_1) \cdots \left [{(\alpha_i^{(n)} - \alpha_i)}\right ]
\cdots (\alpha_i^{(n)} - \alpha_s)},
$$
with the entry in brackets $[\dots]$ excluded.

By assumption $\alpha_i^{(n)}$ converges to $\alpha_i$, and hence also $R_n$ converges
to some non-zero constant $R$.    Given $\delta > 0$, let $N_1 \ge N_0$ be such that
\begin{eqnarray}\label{quantity1}
|\mathrm{Arg}(R) - \mathrm{Arg}(R_n)| < \frac{\delta}{2}
\end{eqnarray}
and
\begin{eqnarray}\label{quantity2}
|\mathrm{Arg}(\alpha_i) - \mathrm{Arg}(\alpha_i^{(n)})| < \frac{\delta}{2n},
\end{eqnarray}
for all $n \ge N_1$.
Then, we have
\begin{eqnarray*}
\mathrm{Arg}(\alpha_i^{(n)} - \alpha_i) &=& \mathrm{Arg}(R_n) - n\mathrm{Arg}(\alpha_i^{(n)})\\
&=& \mathrm{Arg}(R) - n\mathrm{Arg}(\alpha_i) + \delta_n
\end{eqnarray*}
where $\delta_n$ is the sum of the left sides of (\ref{quantity1}) and (\ref{quantity2}).
This proves the claim, with $c = \mathrm{Arg}(R)$.
\qed

\noindent
{\sc Example:}
Let 
$$
P(x) = x^3 + x^2 - 1.
$$
Then $P(x)$ is irreducible and
has exactly two roots $\alpha$ and $\overline\alpha$ outside $C$.
We claim that $\mathrm{Arg}(\alpha)$ is irrational.  Consider the ratio
$$
\omega = {\alpha}/{\overline\alpha}.
$$
Then, since the Galois group of $P(x)$ over the rationals is $S_3$, $\omega$ must
have an algebraic conjugate not on the unit circle, for example,
$$
{\beta}/{\overline\alpha},
$$
where $\beta$ is the real root of $P(x)$.   Thus, $\omega$ is not a root of unity.  
Since $\mathrm{Arg}(\omega) = 2 \mathrm{Arg}(\alpha)$, it follows that
$\mathrm{Arg}(\alpha)$ is irrational.  Thus, by Theorem~\ref{Arg-Theorem}, 
 the relative angle of $\alpha_i^{(n)}$ to $\alpha_i$ is uniformly distributed
as a sequence in $n$.  

Let $\mathrm{Re}(z)$ denote the real part of $z$.  The dot product between two vectors
$\overrightarrow{0 z}$ and $\overrightarrow{0 w}$ is $\mathrm{Re}(z\overline w)$.
It follows from the above that there is no arithmetic progression $kn + \ell$, so that
the sign of
$$
\mathrm{Re}\left [ (\alpha_i^{(kn + \ell)} - \alpha_i) \overline {\alpha_i} \right ]
$$
is constant as a sequence in $n$.  Therefore, $M(Q_n) = \lambda(Q_n)^2$ cannot
be monotone for any arithmetic progression in $n$.

\subsection{Perron polynomials}  

We will show that
for the Salem-Boyd sequence $Q_n(t)$ associated to a Perron polynomial, $\lambda(Q_n)$
is eventually monotone, and prove Theorem~\ref{PerronCase-Theorem}.

\medskip
\noindent
{Proof of Theorem~\ref{PerronCase-Theorem}.}  Let $P$ be a Perron polynomial, and let 
$Q_n(t)$ be an associated
Salem-Boyd sequence.  Let $\mu_1,\dots,\mu_s$ be the roots (counted with
multiplicity) of 
$P$ outside $C$, with $|\mu_1| > |\mu_i|$ for all $i=2,\dots,s$.  By multiplying $P$
by a large enough power of $t$ (this doesn't change $P_*$), we can assume
that $Q_n$ has roots $\lambda^{(n)}_1,\dots,\lambda^{(n)}_s$ outside $C$,
and $|\lambda^{(n)}_i - \mu_i| < | \lambda^{(n)}_i - \mu_j|$ for $\mu_i\neq \mu_j$, and that
$Q_n$ is Perron for all $n \ge 1$.
Let $\lambda^{(n)}_1$ be the largest root of $Q_n$.  Then for all $n$, the root of
$P$ closest to $\lambda^{(n)}_1$ is $\mu_1$, and the root of $Q_n$ closest to $\mu_1$
is $\lambda^{(n)}_1$.   This also implies that $\lambda^{(n)}_1$ is a simple root of $Q_n$.
Fixing $n$, we will show that $\lambda^{(n+1)}_1$ lies strictly between $\lambda^{(n)}_1$
and $\mu_1$.

Consider the equations
\begin{eqnarray}\label{Q-eqn}
0 = Q_{n}(\lambda^{(n)}_1) = (\lambda^{(n)}_1)^{n} P(\lambda^{(n)}_1) \pm P_*(\lambda^{(n)}_1),
\end{eqnarray}
and
$$
Q_{n+1}(\mu_1) = \pm P_*(\mu_1) = Q_n(\mu_1).
$$
Since each of the $Q_n$ are increasing for $t > \lambda^{(n)}_1$, and  $Q_n$
does not have any roots strictly between  $\mu_1$ and  $\lambda^{(n)}_1$, it follows that 
the sign of  $\mu_1 - \lambda^{(n)}_1$ equals the sign of $\pm P_*(\mu_1)$
and does not depend on $n$.

Suppose $\lambda^{(n)}_1 < \mu_1$.  Then, using (\ref{Q-eqn}) in the second
line below, we have
\begin{eqnarray*}
Q_{n+1}(\lambda^{(n)}_1) &=& (\lambda^{(n)}_1)^{n+1} P(\lambda^{(n)}_1) 
\pm P_*(\lambda^{(n)}_1)\\
&=&\lambda^{(n)}_1 (\mp P_*(\lambda^{(n)}_1) \pm P_*(\lambda^{(n)}_1)\\
&=&\pm P_*(\lambda^{(n)}_1)(1 - \lambda^{(n)}_1).
\end{eqnarray*}
By assumption $\lambda^{(n+1)}_1 > 1$.  Also, $P(\lambda^{(n)}_1) < 0$, since
otherwise $P$ would have a real root between $\lambda^{(n)}_1$ and $\mu_1$,
contradicting the assumption that $\lambda^{(n)}_1$ is closer to $\mu_1$
than any other root of $P$.  This implies that
$\pm P_*(\lambda^{(n)}_1) >  0$, and hence $Q_{n+1}(\lambda^{(n)}_1) < 0$, and
$\lambda^{(n)}_1 < \lambda^{(n+1)}_1$.

If $\lambda^{(n)}_1 > \mu_1$, then  $P(\lambda^{(n)}_1) > 0$, and hence 
$\pm P_*(\lambda^{(n)}_1) < 0$.  We thus have
$$
Q_n(\lambda^{(n+1)}_1) = \pm P_*(\lambda^{(n)}_1)(1- \frac{1}{\lambda^{(n+1)}_1}) < 0,
$$
and $\lambda^{(n)}_1 > \lambda^{(n+1)}_1$.\qed

The monotonicity property of Salem-Boyd sequences $Q_n$ associated to
a Perron polynomial $P$ allows us to give
a lower bound greater than one for the sequences $\lambda(Q_n)$.

\begin{proposition}\label{lower-bound-proposition}  If $Q_n(t)$ is defined by 
$$
Q_n(t) = t^n P(t) \pm P_*(t),
$$
where $P$ is a Perron polynomial, and $n_0$ is such that
$\lambda(Q_n)$ is monotone for $n \ge n_0$, then
$$
\lambda(Q_n) \ge \min\{\lambda(Q_{n^\pm_0}),\lambda(P)\}
$$
for all $n \ge n_0$. \end{proposition}

\subsection{P-V  and Salem polynomials.}\label{PVcase-subsection}

We now consider the case when $P=P^\pm_{\Sigma,\tau}$ belongs to a special
class of Perron polynomials, 
namely those satisfying  $N(P^\pm_{\Sigma,\tau}) = 1$.  

 A {\it P-V number} is a real algebraic integer
$\alpha > 1$ such that all other algebraic conjugates lie strictly within $C$.
  A {\it Salem number} is a real algebraic integer $\alpha > 1$ such that
all other algebraic conjugates lie on or within $C$ with at least one on $C$.
If $f$ is an irreducible monic integer polynomial with $N(f) = 1$, then the root of $f$ outside
$C$ has absolute value equal to either a Salem number, if $f$ has degree greater
than 2 and is reciprocal, or a P-V number otherwise.   
 If $f$ is reciprocal and $N(f) = 1$l, then $\lambda(f)$ is either a Salem number
or a quadratic P-V number.

The polynomials $Q^\pm_n(t)$ were originally studied by Salem \cite{Salem44} in the 
case when $P(t)$ is a P-V polynomial to show that every P-V number is the upper
and lower limit of Salem numbers.  Boyd \cite{Boyd77} showed that any Salem number
occurs as $M(Q^\pm_n)$ for some P-V polynomial $P(t)$.    

Assume that $P(t)$ has no reciprocal factors and $P(1) \neq 0$.  Let 
$$
n^-_0(P) = d - 2\frac{P'(1)}{P(1)} + 1
$$
where $d$ is the degree of $P$,
and let
$$
n^+_0(P) = 1
$$
for all $P$.   
For any polynomial  (or Laurent polynomial) $P$, let $\ell(P)$ be the sign of the lowest degree
coefficient of $P$.  The following Proposition is proved in Boyd's discussion in 
(\cite{Boyd77} p. 320-321), and implies Theorem~\ref{SalemCase-Theorem}.

\begin{proposition}\label{PVcase-proposition} If $P$ is a P-V polynomial for the P-V number 
$\theta$,
then the polynomial $Q^\pm_n(t)$ has a real root greater than
one if and only if $n \ge n^\pm_0(P)$.   Furthermore, the sequences of resulting Salem numbers
$\alpha^\pm_n$ is monotone increasing (decreasing) if and only if
$\pm \ell(P) > 0$ ($<0$).
\end{proposition}

\noindent{Proof.}  The proof follows from looking at the real graphs of $Q^\pm_n(t)$
and of $P$.  Since $P(1) = P_*(1) < 0$, $Q^+_n(1)$ must be strictly negative.   Thus,
$Q^+_n$ must have a root larger than 1 for all $n$, and we can set $n^+_0 = 1$.
The graph of $y=Q^-_n(t)$ passes through the real axis at $t=1$.  Thus, $Q^-_n(t)$
has a positive real root if and only if the derivative of $Q^-_n$ is negative.  Note that
$Q^-_n(t)$ cannot have a negative real root by the argument in the proof of
Lemma~\ref{numberofroots-lemma}.  This proves
the first part of the Proposition.  

For the second part, note that since $P$ has only 
one root $\theta$ outside $C$, $P_*(\theta)$ and $\pm \ell(P)$ must have the same 
sign.  Suppose, for example, that $\pm\ell(P) > 0$.  Put $\alpha^\pm_n = \lambda(Q^\pm_n)$.
Then $Q^\pm_n(\theta) > 0$, and 
hence $\theta > \alpha^\pm_n$ for all $n$.  This implies that $P(\alpha_{n+1}^\pm) < 0$.
Now consider the equations:
\begin{eqnarray*}
Q^\pm_n(\alpha^\pm_{n+1}) &=& Q^\pm_n(\alpha^\pm_{n+1}) - Q^\pm_{n+1}(\alpha^\pm_{n+1})\\
&=& \left ((\alpha^\pm_{n+1})^{n} -( \alpha^\pm_{n+1})^{n+1} \right ) P(\alpha^\pm_{n+1}).
\end{eqnarray*}
The bottom formula is a product of negative numbers.  Hence, $Q^\pm_n(\alpha^\pm_n) > 0$, and
$\alpha^\pm_{n+1} > \alpha^\pm_n$.\  The case $\pm\ell(P)<0$ is proved in an analogous way.\qed

\section{Poset structure on fibered links}\label{applications-section}

We now apply results of the previous sections to sequences of fibered links
obtained by iterated trefoil plumbings.   Let $(K,\Sigma)$ be a fibered link,
and let $P$ be the polynomial produced by a given locus of plumbing $\tau$. 
Let $\Delta_n = \Delta_{(K_n,\Sigma_n)}$ be the Alexander polynomials
of the iterated trefoil plumbings.  If $P$ is a Perron polynomial, 
then Proposition~\ref{lower-bound-proposition} implies that one can find lower
bounds for $\lambda(\Delta_n)$, and hence for $M(\Delta_n)$ at least for large $n$.  The 
situation is even better when $P$ is a P-V polynomial.  In this case, we can explicitly 
find the minimal $\lambda(\Delta_n)$ and hence $M(\Delta_n)$
 in the sequence by comparing $\lambda(\Delta_{n_0})$
and $\lambda(P)$, where $n_0$ is as in Proposition~\ref{PVcase-proposition}.
Furthermore, any P-V polynomial  satisfies the
inequality  (see \cite{Siegel44})
$$
\lambda(P) \ge \lambda(x^3-x-1) \approx 1.32472.
$$
It is not known in general if there is a lower bound greater than one for Salem numbers.

A fibered link $(K,\Sigma)$ will be called a {\it Salem fibered link}, if the following equivalent statements hold:
\begin{description}
\item{(1)} $N(\Delta_{(K,\Sigma)}) = 1$;
\item{(2)} $\lambda(\Delta_{(K,\Sigma)}) = M(\Delta_{(K,\Sigma)})$; and
\item{(3)} $M(\Delta_{(K,\Sigma)})$ is a Salem number or a quadratic P-V number.
\end{description}

Let $\sS$ be the set of Salem fibered links, and write
$$
(K_1,\Sigma_1) \prec_S (K_2,\Sigma_2)
$$
if $(K_2,\Sigma_2)$ can be obtained from $(K_1,\Sigma_1)$ be a sequence of trefoil
plumbings, where the polynomial $P^\pm_{\Sigma,\tau}$ corresponding to the plumbing
locus at each stage is a P-V polynomial.  If $(K_1,\Sigma_1) \prec_S (K_2,\Sigma_2)$,
then the topological Euler characteristic of $\Sigma_1$ is strictly less than that of
$\Sigma_2$.  Thus, $\prec_S$ defines an (anti-symmetric) partial order on Salem 
fibered links.  Proposition~\ref{lower-bound-proposition}
implies the following.

\begin{proposition}\label{poset-prop} If $(K_1,\Sigma_1) \prec_S (K_2,\Sigma_2)$, then
$$
M(\Delta_{(K_2,\Sigma_2)}) \ge \min \{ M(\Delta_{(K_1,\Sigma_1)}), \theta_0\}
$$
where $\theta_0 \approx 1.32472$ is the smallest P-V number.
\end{proposition}

Consider the graph structure of $\sS$ with respect $\prec_S$.    By Proposition~\ref{poset-prop},
for any connected 
subgraph of $\sS$, the minimal Salem number
can be determined by comparing the minimal elements with respect to $\prec_S$.

\begin{question} Is $\sS \cap \sK$ connected with respect to $\prec_S$?  
\end{question}

It is not difficult to produce examples of Salem fibered links $(K,\Sigma)$ and a locus for plumbing
$\tau$ 
such that $P_{\Sigma,\tau}$ is not a P-V polynomial (see Section~\ref{example-section}).     
We will say a Salem fibered link 
$(K,\Sigma) \in \sS \cap \sK$ is {\it isolated}  if for all loci of plumbing $\tau$ on $\Sigma$,
 the corresponding polynomial $P$ is not a P-V polynomial.
  
\begin{question} Are there isolated Salem links?  
\end{question}

Although we do not know of any isolated Salem links, Salem fibered links do appear
sporadically in Salem-Boyd sequences not associated to P-V polynomials
as seen in the table at the end of Section~\ref{example-section}.

\section{A family of fibered two bridge links.}\label{example-section}

\begin{figure}[htbp]
\begin{center}
\includegraphics[height=1.5in]{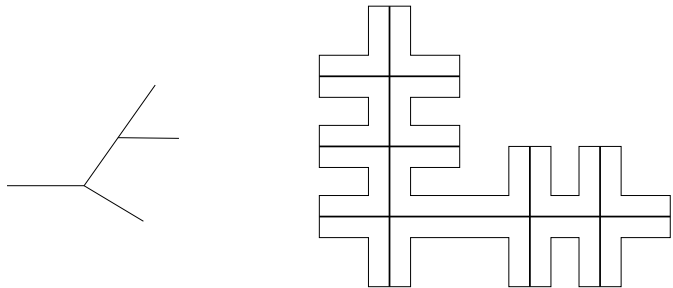}
\caption{Construction of fibering surface for arborescent link}
\label{surfaceplumbing-fig}
\end{center}
\end{figure}

The simplest examples to consider are those coming from arborescent
links.  Let $\Gamma$ be a tree, with vertices $\nu$ with  labels $m(\nu) = \pm 1$.
Let $\sL$ be a
union of line segments in the plane, intersecting transversally, whose dual graph is $\Gamma$,
and let $U(\sL)$ be the surface obtained by thickening $\sL$.  This is illustrated in
Figure~\ref{surfaceplumbing-fig}.

 Consider the surface in Figure~\ref{surfaceplumbing-fig}  as a subspace of $S^3$ and glue 
 together opposite sides in the diagram that are connected by a vertical or horizontal path
 with a positive or negative
 full-twist according to the labeling on the graph.  The resulting surface $\Sigma$ is a fibering
 surface for $K = \partial \Sigma$ by \cite{Stallings:fibered}, since it can be obtained by
a sequence of Hopf plumbings on the unknot.   The line segments of $\sL$ close up to form
a free basis for $\HH_1(\Sigma;\R)$.    Thus, the vertices of $\Gamma$ can be thought of
as basis elements of $\HH_1(\Sigma;\R)$.  Let $S_\Gamma$ be the matrix where the rows
and columns correspond to vertices $\nu_1,\dots,\nu_k$
 of $\Gamma$, and the entries $a_{i,j}$ are given by
 $$
 a_{i,j} = 
 \left \{
 \begin{array}{rl}
 -1 &\qquad\mbox{if  $i< j$, and
 $\nu_i$ and $\nu_j$ are connected by an edge}\\
m(\nu_i) &\qquad\mbox{if $i=j$, and}\\
0 &\qquad \mbox{otherwise.}
\end{array}
\right .
$$
Then  $S$ is a Seifert matrix for $(K,\Sigma)$.  It follows that although there may be
several fibered links $(K,\Sigma)$ associated to a given labeled graph $\Gamma$, the Seifert matrix,
and hence the Alexander polynomial, is determined by $\Gamma$.

\begin{figure}[htbp]
\begin{center}
\includegraphics[width=3in]{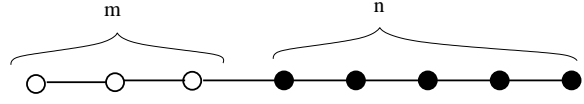}
\caption{Plumbing graph with positive (negative) vertices filled black (white)}
\label{bamboo-fig}
\end{center}
\end{figure}

\begin{figure}[htbp]
\begin{center}
\includegraphics[height=1.5in]{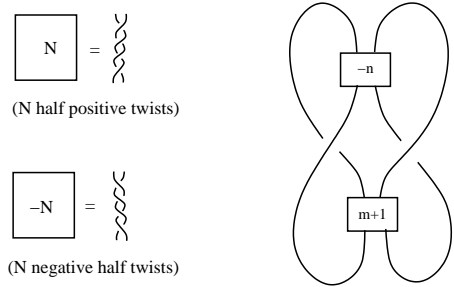}
\caption{Two bridge link associated to $\Gamma_{m,n}$}
\label{twobridge-fig}
\end{center}
\end{figure}

Consider the family of examples $\Gamma_{m,n}$ 
in Figure~\ref{bamboo-fig}.
The associated fibered links $(K_{m,n},\Sigma_{m,n})$ (determined uniquely by
$\Gamma_{m,n}$) are 
the two-bridge link drawn in Figure~\ref{twobridge-fig}.

Fixing $m$, and letting $n$ vary gives a sequence of fibered links $(K_{m,n},\Sigma_{m,n})$ that
are obtained by iterated plumbing on $(K_{m,1},\Sigma_{m,1})$.   Thus, the Alexander polynomials
$\Delta_{m,n} = \Delta_{K_{m,n},\Sigma_{m,n}}$ are Salem-Boyd sequences associated to some polynomials
$P_m$.  We will compute the $P_m$, and their numerical invariants.

Considering the vertices of $\Gamma_{m,1}$ as basis elements in $\HH_1(\Sigma_{m,1},\R)$,
the path $\tau$ is dual to the right-most vertex.  We start with $\Gamma_{1,1}$.
The link $K_{1,1}$  is the figure-eight knot, or 4${}_1$ in 
Rolfsen's table \cite{Rolfsen76}.  
 We will use Equation~\ref{FirstP-equation} to find $P_1$.
Thus, 
$P_1$ is given by
\begin{eqnarray*}
P_1 (t) &=& \sign(S) \left | 
t\ \left ( 
\begin{array}{cc}
-1 & 0\\
-1 & 1
\end{array}
\right )
-
\left (
\begin{array}{cc}
-1 & -1 \\
0 & 0
\end{array}
\right )\right |\\
&=& t(t-2)
\end{eqnarray*}
Since $P_1$ has only one root outside $C$, we have the following Proposition.

\begin{proposition} The links $(K_{1,n},\Sigma_{1,n})$ are Salem fibered links.
\end{proposition}

The Salem numbers $\lambda(\Delta_{1,n})$ converge to $\lambda(P) = 2$, from above
for $n$ odd, and from below for $n$ even.   The smallest Salem number in this sequence
occurs for $(K_{1,4},\Sigma_{1,4})$, and is approximately 1.8832.   

From $P_1$ it is possible to compute all the $P_m$ using Equation~\ref{SecondP-equation}.
We first recall that $(K_{m,0},\Sigma_{m,0})$ is the $(2,m+1)$ torus link, $T_{(2,m+1)}$.
The Alexander polynomial is given by
$$
\Delta_{m,0} (t) = \frac{t^{n+1} + (-1)^n}{t+1}.
$$
Since $P_1(t) = t(t-2)$, and $K_{1,1}$ has one component, we also have
\begin{eqnarray*}
\Delta_{1,n}(t) &=& \frac{ t^n P_1(t) + (-1)^{n+1}(P_1)_*(t)}{t+1}\\
&=& \frac{t^{n+1}(t-2) + (-1)^{n+1}(-2t+1)}{t+1}\\
&=& \frac{t^{n+1}(t-2) + (-1)^n 2t + (-1)^{n+1}}{t+1}
\end{eqnarray*}

Furthermore, $\Gamma_{m,0}$ can be thought of as a subgraph of $\Gamma_{m,1}$,
and if $S_{m,0}$ and $S_{m,1}$ are their associated Seifert surfaces, we have
$$
\sign(S_{m,0}) = \sign(S_{m,1}).
$$
By Equation~\ref{SecondP-equation}, we have
\begin{eqnarray*}
P_m (t) &=& \Delta_{m,1}(t)  + \Delta_{m,0}(t)\\
&=&\frac{t^{m+1}(t-2) + (-1)^m 2t + (-1)^{m+1} + t^{m+1} + (-1)^m}{t+1}\\
&=&\frac{t^{m+2} - t^{m+1} + (-1)^m 2t}{t+1}\\
&=&\frac{t(t^m(t-1) + (-1)^m 2)}{t+1}
\end{eqnarray*}
Since we are only concerned with $\Delta_{m,n}$ and hence
$P_m$ up to products of cyclotomic polynomials, it is
convenient to rewrite $P_m$ as 
$$
P_m(t) = t(t^m(t-1) + (-1)^m 2).
$$

\begin{proposition}\label{rootsofP-prop}
All roots of $P_m(t)$ other than $0$ and $-1$ lie outside $C$, hence 
$$
M(P_m) = 2\qquad{\mathrm{and}}\qquad N(P_m) = m.
$$
\end{proposition}

\medskip
\noindent{Proof.}
Suppose $|t| \leq 1$, then $|t^m(t-1)| \leq 2$ with equality if and only if $t = -1$.\qed

\begin{proposition}\label{radofP-prop} 
$$
\lim_{m \rightarrow \infty} \lambda(P_m) = 1.
$$
\end{proposition}

\medskip
\noindent{Proof.} Take any $\epsilon > 0$.  Let $D_\epsilon = \{ z \in \C \ : \ |z| > 1+\epsilon\}$.
Let $\overline{D_\epsilon}$ be the closure of $\C$  in the Riemann sphere.
Then for large $m$ 
$$
\frac{2}{|t^m|} < \frac{|t-1|}{|t|}
$$
for all $t$ on the boundary of $\overline{D_\epsilon}$ and both sides are analytic
on $\overline{D_\epsilon}$.  Therefore, by Rouch\'e's theorem $P_m$ has no 
roots on $\overline{D_\epsilon}$ for large $m$.\qed

\begin{corollary} The homological dilatations of $(K_{m,n},\Sigma_{m,n})$ can be
made arbitrarily small by taking $m$ and $n$ large enough.
\end{corollary}

Salem fibered links appear sporadically as homological
dilatations of $(K_{m,n},\Sigma_{m,n})$
for $m,n > 1$.  
A list for $1<m,n< 60$ found by computer search is given in the table below.  
The minimal polynomials, which are reciprocal, are denoted by a list of the first half of
the coefficients.

\begin{center}
\begin{tabular}{|c|l|l|}
\hline
$(m,n)$ & Salem number & Minimal polynomial\\
\hline
(3,5) & 1.63557 & 1\ -2\ 2\ -3\\
\hline
(3,8) & 1.50614 & 1\ -1\ 0\ -1\\
\hline 
(5,9) & 1.42501 &1\ -1\ 0\ -1\ 1\\
\hline
\end{tabular}
\end{center}

\begin{question}
Are the Salem fibered links in the table above isolated in the sense 
of Section~\ref{applications-section}?
\end{question}

Salem numbers also appear as roots of irreducible 
factors of the Alexander polynomial.  For example,
the Alexander polynomial for $K_{11,21}$ has largest root equal to the 7th
smallest known Salem number \cite{Moss-url}.   Its minimal polynomial is
given by 
$$
\Delta_{K_{11,21}}(x) = x^{10} - x^7 - x^5 - x^3 + 1.
$$
The monodromy $h_{m,n}$ of the fibered links $(K_{m,n},\Sigma_{m,n})$
were also studied by
Brinkmann \cite{Brinkmann04}, who showed that $h_{m,n}$
is  pseudo-Anosov for all $m,n$, and that the dilatations
converge to 1 as $m,n$ approach infinity.

 \bibliographystyle{plain}
 \bibliography{math}
\bigskip

\verse{
Eriko Hironaka\\
Department of Mathematics\\
Florida State University\\
Tallahassee, FL 32306-4510\\
U.S.A.
}

  \end{document}